\documentclass[a4paper]{amsart}
\usepackage[latin1]{inputenc}
\usepackage[english]{babel}

\usepackage{emptypage}

\usepackage{amsmath}
\usepackage{amssymb}
\usepackage{amsfonts}
\usepackage{mathrsfs}
\usepackage{amsthm}
\usepackage{doi}

\usepackage[final]{showkeys}

\usepackage{color}
\usepackage{graphicx}

\usepackage[all]{xy}

\usepackage{cite}
\usepackage{epigraph}

\usepackage{geometry}

\makeatletter
\newtheorem*{rep@theorem}{\rep@title}
\newcommand{\newreptheorem}[2]{%
\newenvironment{rep#1}[1]{%
 \def\rep@title{#2 \ref{##1}}%
 \begin{rep@theorem}}%
 {\end{rep@theorem}}}
\makeatother

\newcommand{\C}{\mathbb{C}}
\newcommand{\rg}{\mathrm{reg}}

\newcommand{\supp}{\mathrm{supp}\,}

\renewcommand{\H}{\mathcal{H}}
\newcommand{\R}{\mathbb{R}}
\newcommand{\Z}{\mathbb{Z}}
\newcommand{\Ci}{\mathcal{C}}
\newcommand{\E}{\mathcal{E}}
\newcommand{\Di}{M}
\newcommand{\D}{M}
\newcommand{\Dii}{\mathscr{D}}
\newcommand{\B}{\mathscr{B}}
\newcommand{\N}{\mathbb{N}}
\newcommand{\EL}{\mathfrak{L}}

\newcommand{\Ol}{\mathcal{O}}
\newcommand{\DD}{\mathbb{D}}

\newcommand{\Lip}{\mathrm{Lip}}

\newcommand{\mass}{\mathbf{M}}

\newcommand{\de}{\partial}
\newcommand{\debar}{\overline{\de}}

\theoremstyle{plain}
\newtheorem{Teo}{Theorem}[section]
\newreptheorem{Teo}{Theorem}
\newtheorem{Prp}[Teo]{Proposition}
\newtheorem{Cor}[Teo]{Corollary}
\newtheorem{Lmm}[Teo]{Lemma}
\theoremstyle{definition}
\newtheorem{Def}{Definition}[section]
\theoremstyle{remark}
\newtheorem{Rem}{Remark}[section]

\title{Positive metric currents and holomorphic chains in Hilbert spaces}
\author{Samuele Mongodi}
\address{Scuola Normale Superiore di Pisa - Piazza dei Cavalieri 7, I-56126 Pisa}
\email{s.mongodi@sns.it}

\begin{document}

\begin{abstract}We present some results concerning currents of integration on finite-dimensional analytic spaces in Hilbert spaces, using the setting of metric currents. In particular, we obtain the characterization of such currents as positive closed $(k,k)-$rectifiable currents and solve the boundary problem for holomorphic chains.
\end{abstract}

\maketitle

\section{Introduction}

The theory of currents found deep and important applications in complex analysis and geometry; just to mention a few of them, we recall the characterization of holomorphic chains by King in \cite{king1} and by Harvey and Shiffman in \cite{harvey2}, the removal of singularities for analytic functions and sets by Shiffman in \cite{shiffman1} and the boundary problem for holomorphic chains by Harvey and Lawson in \cite{harvey1, harvey3}.

\medskip

Such a powerful tool was extended to general metric spaces in \cite{ambrosio1} and we described some possible applications to complex analysis in \cite{mongodi1}, specializing the theory in the case of singular complex spaces and complex Banach spaces. In order to gain a wider understanding of the latter, here we turn our attention to holomorphic metric chains in Hilbert spaces, tackling two specific problems: firstly, we consider positive $(k,k)$-rectifiable closed metric currents and investigate their link with complex analytic sets; then, we proceed to study the boundary problem for holomorphic (metric) chains and give a positive answer under suitable hypotheses of regularity and of general position.

\medskip

We organized the material as follows. 

After recalling the basic definitions and properties of metric currents and rectifiable sets in metric spaces (from \cite{ambrosio1, ambrosio2}), we summarize the results presented in \cite{mongodi1} on the extension of metric currents to complex Banach spaces.

In Section 4, the properties of analytic sets and of the currents of integration associated to them are investigated; in particular, we introduce the concept of positive metric currents, presenting some simple properties, and we obtain a weak analogue of Wirtinger's formula.

The main result of this section is the extension of King's characterization theorem (see \cite{king1}) for holomorphic chains to Hilbert spaces.
\begin{repTeo}{teo_king}Let $\Omega\subset H$ be a ball, $S$ be a rectifiable current in $\Omega$. Suppose that
\begin{enumerate}
\item $\supp dS\cap\Omega=\emptyset$;
\item $S$ is a $(k,k)$ positive current.
\end{enumerate}
Then $S$ can be represented as a sum with integer coefficients of integrations on the regular parts of analytic sets. \end{repTeo}

In the last section, we turn our attention to  maximally complex metric currents and CR-manifolds. We recover a finite-dimensional embedding result and a characterization for maximally complex integration currents.

Finally, we tackle the boundary problem for holomorphic chains in a Hilbert space: we solve the problem with a technical hypothesis on $M$, specifically, the existence of a finite rank projection with transverse self-intersections.

\begin{repTeo}{teo_HL1}Let $M$ be a compact, oriented $(2p-1)-$manifold (without boundary) of class $\Ci^2$ embedded in $H$, and suppose that there exists an orthogonal decomposition $H=\C^p\oplus H'$ such that the projection $\pi:H\to\C^p$, when restricted to $M$, is an immersion with transverse self-intersections. Then if $M$ is an $MC-$cycle there exists a unique holomorphic $p-$chain $T$ in $H\setminus M$ with $\supp T\Subset H$ and finite mass, such that $dT=[M]$ in $H$.\end{repTeo}

It is worth noting that the Hilbert space hypothesis isn't clearly required in these results, but our methods rely heavily on the possibility of estimating norms in terms of coordinates, which is a feature essentially related to the existence of an orthonormal basis.

\medskip

\noindent{\bf Acknowledgements. }  I acknowledge the support of the Italian project FIRB-IDEAS ''Analysis and Beyond''. I also thank prof. Luigi Ambrosio for the fruitful discussions about the geometric measure theory tools employed in the proof of King's result.

\section{Metric currents}

Let $X$ be a metric space. Let us denote by $\Lip(X)$ the space of complex-valued Lipschitz functions on $X$ and by $\Lip_b(X)$ the algebra of bounded complex-valued Lipschitz functions. Following \cite{mongodi1}, we introduce the spaces
$$\E^k(X)=\{(f,\pi_1,\ldots, \pi_k)\ \vert\ f\in\Lip_b(X),\ \pi_j\in\Lip(X),\ j=1,\ldots, k\}\;,\quad \E^0(X)=\Lip_b(X)\;.$$
The elements of these spaces are called \emph{metric forms}.

\begin{Def}A \emph{$k-$dimensional metric current} is a functional $T:\E^k(X)\to\C$ satisfying the following
\begin{enumerate}
\item $T$ is multilinear
\item whenever $(f,\pi^i)\to (f,\pi)$ pointwise, with uniformly bounded Lipschitz constants, then $T(f^i,\pi^i)\to T(f,\pi)$
\item $T(f,\pi)=0$, whenever there is an index $j$ for which $\pi_j$ is constant on a neighborhood of $\supp f$
\item $T$ has finite mass, i.e. there exists a finite Radon measure $\mu$ on $X$ such that
$$|T(f,\pi)|\leq \prod_{j=1}^k\Lip(\pi_j)\int_{X}|f|d\mu\qquad \forall (f,\pi)\in\E^k(X)\;;$$
the infimum of such measures $\mu$ is called the \emph{mass} (measure) of $T$.
\end{enumerate}
\end{Def}

We will denote by $M_k(X)$ the space of $k-$dimensional metric currents on $X$; endowing $M_k(X)$ with the mass norm $\|T\|=\mu(X)$, where $\mu$ is the mass measure of $T$, we turn it into a Banach space. We will sometimes write $T(fd\pi)$ for $T(f,\pi_1,\ldots, \pi_k)$.

We refer to \cite{ambrosio1} for a detailed discussion of the general properties of metric currents; here we recall only the main definitions and facts.

\begin{Def}We say that $T_h\in M_k(X)$ \emph{converge weakly} to $T\in M_k(X)$ if
$$T_h(fd\pi)\to T(fd\pi)\qquad \forall\ fd\pi\in \E^k(X)\;.$$
\end{Def}

\begin{Def} Given $T\in M_k(X)$, $k\geq1$, we can define the \emph{boundary} of $T$ by 
$$dT(f,\pi_1,\ldots, \pi_{k-1})=T(1,f,\pi_1,\ldots, \pi_k)\qquad \forall\ fd\pi\in \E^{k-1}(X)\;.$$\end{Def}
The functional $dT$ satisfies the first three assumptions, but will not in general be of finite mass; in that case, we say that $T$ is a \emph{normal} metric current. The space $N_k(X)$ of $k-$dimensional normal metric currents is a Banach space when endowed with the norm $\|T\|_N=\|T\|+\|dT\|$.

\begin{Def}Given a Lipschitz map between complete metric spaces $F:X\to Y$ and $T\in M_k(X)$, we define the \emph{pushforward} $F_\sharp T\in M_k(Y)$ of $T$ through $F$ by
$$F_\sharp T(f,\pi_1,\ldots,\pi_k)=T(f\circ F, \pi_1\circ F, \ldots, \pi_k\circ F)\qquad \forall (f,\pi)\in \E^k(Y)\;.$$\end{Def}

We observe that, if $\mu$ and $\widetilde{\mu}$ are the mass measures of $T$ and $F_\sharp T$ respectively, we have
$$\widetilde{\mu}\leq \Lip(F)^{k}F_*\mu\;.$$
Moreover, one has $F_\sharp(dT)=d(F_\sharp T)$.

\begin{Def}Given $T\in M_k(X)$ and $\omega=(u,v_1,\ldots, v_h)\in \E^h(X)$, with $h\leq k$, the \emph{contraction} of $T$ with $\omega$ is $T\llcorner\omega\in M_{k-h}(X)$ and it is defined by
$$T\llcorner\omega(f,\pi_1,\ldots, \pi_{k-h})=T(fu, v_1,\ldots, v_h, \pi_1,\ldots, \pi_{k-h})$$
for every $(f,\pi_1,\ldots, \pi_{k-h})\in \E^{k-h}(X)$.\end{Def}

We have that
$$\|T\llcorner \omega\|\leq \sup|u|\prod_{j=1}^h\Lip(v_j)\|T\|\;.$$

\begin{Def}The \emph{support} of $T\in M_k(X)$ is the least closed set $\supp T$ such that if $\supp f$ doesn't intersect $\supp T$, then $T(f,\pi)=0$.\end{Def} 

By the finiteness of the mass, we can extend a metric current as a functional on $\B^\infty(X)\times[\Lip(X)]^k$, that is, we can allow the first entry of the metric forms to be a bounded Borel function on $X$. The previous definitions and remarks all extend without changes.

\begin{Prp}Let $T\in M_k(X)$ be a metric current, then
\begin{enumerate}
\item $T$ is alternating in $\pi_1,\ldots, \pi_k$,
\item $T$ satisfied the chain rule and the Leibniz rule,
\item $T=i_\sharp S$ with $S\in M_k(\supp T)$ and $i:\supp T\to X$ the inclusion.
\end{enumerate}
\end{Prp}

The proofs of these statements can be found in \cite{ambrosio1}.

We recall a compactness result for normal currents.

\begin{Teo}Let $(T_h)\subset N_k(X)$ be a bounded sequence and assume that for any integer $p\geq 1$ there exists a compact $K_p\subset E$ such that
$$\|T_h\|(E\setminus K_p)+\|dT_h\|(E\setminus K_p)<\frac1p\qquad\forall\ h\in \N\;.$$
Then, there exists a subsequence $(T_{h(n)})$ converging to a current $T\in N_k(X)$ such that
$$\|T\|(E\setminus \bigcup_{p=1}^\infty K_p)+\|dT\|(E\setminus\bigcup_{p=1}^\infty K_p)=0\;.$$\end{Teo}

Finally, we recall the comparison theorem between classical and metric currents, in the form given in \cite{mongodi1}.

\begin{Teo}\label{teo_comp_man}Let $U$ be an $N-$dimensional complex manifold, $N\geq 1$, endowed with the distance given by an hermitian metric. For every $m\geq0$ there exists an injective linear map $C_m:\Di_m(U)\to \Dii_m(U)$ such that
$$C_m(T)(fdg_1\wedge\ldots\wedge dg_m)=T(f,g_1,\ldots, g_m)$$
for all $(f,g_1,\ldots, g_m)\in\Ci^\infty_c(U)\times[\Ci^\infty(U)]^m$. The following properties hold:
\begin{enumerate}
\item for $m\geq1$, $d\circ C_m=C_{m-1}\circ d$;
\item there exists a positive constant $c_1$ such that, for all $T\in\Di_m(U)$, 
$$c_1^{-2}\|T\|\leq\mass(C_m(T))\leq c_1^{2}{N\choose m}\|T\|\;$$
\item the restriction of $C_m$ to $N_{m}$ is an isomorphism onto $\mathfrak{N}_{m}$;
\item the image of $C_m$ contains the space $\mathfrak{F}_{m}(U)$.\end{enumerate}
\end{Teo}

\subsection{Rectifiable currents and slicing}

Let $\H^k$ be the $k-$dimensional Hausdorff measure on $X$.
\begin{Def}\label{def_rect_set}We say that a $\H^k-$measurable set $S\subset X$ is \emph{countably $\H^{k}-$rectifiable} if there exist sets $A_i\subset\R^k$ and Lipschitz functions $f_i:A_i\to X$ such that 
$$\H^k\left(S\setminus\bigcup_{i=1}^\infty f_i(A_i)\right)=0\;.$$\end{Def}

\begin{Lmm}Let $S\subset X$ be countably $\H^k-$rectifiable. Then, there exist finitely or countably many  compact sets $K_i\subset\R^k$ and bi-Lipschitz maps $f_i:K_i\to S$ such that their images are pairwise disjoint and $\H^k(S\setminus\cup_if_i(K_i))=0$.\end{Lmm}

\begin{Def}A current $T\in M_{k}(X)$ is said to be \emph{rectifiable} if 
\begin{enumerate}
\item $\|T\|$ is concentrated on a countably $\H^k-$rectifiable set;
\item $\|T\|$ vanishes on $\H^k-$negligible Borel sets.
\end{enumerate}
The space of  rectifiable currents is denoted by $\mathcal{R}_{k}(X)$.

We say that a rectifiable current $T$ is \emph{integer rectifiable} if for any $\phi\in\Lip(X,\R^k)$ and open set  $B$ in $X$ one has $\phi_\sharp(T\llcorner B)=[u]$ with $u\in L^1(\R^k, \Z)$. The space of such currents is denoted by $\mathcal{I}_{k}(X)$.

We define the spaces of \emph{integral} currents as follows.
$$I_k(X)=\mathcal{I}_{k}\cap N_k(X)\;:$$\end{Def}

\begin{Teo}Let $T\in M_{k}(X)$, $k\geq1$. Then $T\in\mathcal{R}_{k}(X)$ (resp. $T\in\mathcal{I}_{k}(X)$) if and only if there exist a sequence $\{K_i\}$ of compact sets in $\R^k$, a sequence $\{\theta_i\}$ of functions in $L^1(\R^k,\R)$ (resp. $L^1(\R^k,\Z)$) with $\supp \theta_i\subset K_i$ and a sequence $\{f_i\}$ of bi-Lipschitz maps $f_i:K_i\to X$ such that 
$$\|T\|(A)=\sum_i\|(f_i)_\sharp[\theta_i]\|(A)\qquad T(f,\pi)=\sum_i(f_i)_\sharp[\theta_i](f,\pi)$$
for every Borel set $A\subset X$ and for every $(f,\pi)\in\E^k(X)$.\end{Teo}

\begin{Def}Given $T\in \mathcal{R}_k(X)$, we set
$$S_T=\{x\in X\ :\ \Theta_k(\|T\|,x)>0\}$$
and we define the \emph{size} of $T$ as
$$\mathbf{S}(T)=\H^k(S_T)\;.$$
\end{Def}
One can show that $S_T$ is countably $\H^k-$rectifiable, that $\|T\|$ is concentrated on $S_T$ and that any other Borel set on which $\|T\|$ is concentrated includes $S_T$ up to $\H^k-$negligible sets.

\medskip

Given $T\in M_{k}(X)$, $\pi\in\Lip(X,\R^m)$, we define the \emph{slice} $\langle T,\pi,x\rangle\in \Di_{k-m}(X)$ by
$$\langle T,\pi,x\rangle(f,\eta)=\lim_{\epsilon\to0}T(f\rho_\epsilon\circ\pi, \pi, \eta)$$
where $\rho_\epsilon$ is any family of mollifiers, for every $x\in\R^m$ for which the limit exists.

\begin{Teo}\label{teo_ext_slice}If $T\in N_{k}(X)$, $\pi\in Lip(X,\R^m)$, then
\begin{enumerate}
\item for $\mathcal{L}^m-$almost every $x\in\R^m$, the slice $\langle T,\pi,x\rangle$ exists and is  normal and $d\langle T,\pi, x\rangle=(-1)^m\langle dT, \pi, x\rangle$;
\item for all $(f,g)\in\mathcal{B}^\infty\times[\Lip(X)]^{k-m}$, 
$$\int_{\R^m}\langle T,\pi,x\rangle(f,g)dx=T\llcorner(1,\pi)(f,g)\;;$$
\item for every $\|T\llcorner(1,\pi)\|-$measurable set $B\subset X$, 
$$\int_{\R^m}\|\langle T,\pi, x\|(B)dx=\|T\llcorner(1,\pi)\|(B)\;.$$
\end{enumerate}
\end{Teo}

\subsection{Rectifiable currents in Banach spaces}

Let us specialize the theory of metric currents to the case when $X$ is a Banach space $(E, \|\cdot\|)$. We recall some notions on metric differentiability, Jacobians, area and coarea formulas from \cite{ambrosio1} and \cite{ambrosio2}.

Let us suppose that $E$ is a $w^*-$separable dual space, i.e. that $E=G^*$ and that there exists a sequence $(\phi_h)_h\subset \Lip_1(G)$ such that
$$\|x-y\|_G=\sup_h|\phi_h(x)-\phi_h(y)|\qquad\forall\ x,y\in G\;.$$
Given instead $x,y\in E$, we define
$$d_w(x,y)=\sum_{n=0}^\infty 2^n|\langle x-y, g_n\rangle|\;,$$
where $\{g_n\}_n$ is a dense subset of the unit ball of $G$. The topology induced by $d_w$ on the bounded sets of $E$ is called $w^*-$topology.

\begin{Def}A sequence $(T_h)\subset M_k(E)$ is said to \emph{$w^*-$converge} to $T\in M_k(E)$ if
$T_h(fd\pi)$ tends to $T(fd\pi)$ for every $fd\pi\in \E^k(E)$ with $w^*-$continuous coefficients.\end{Def}

We have the following weak$\!\!\phantom{|}^*$-compactness result.
\begin{Teo}Let $Y$ be a $w^*-$separable dual space, let $(T_h)\subset N_k(Y)$ be a bounded sequence, and assume that for any $\epsilon>0$ there exists $R>0$ such that $K_h=\overline{B}_R(0)\cap \supp T_h$ are equi-compact and
$$\sup_{h\in \N}\|T_h\|(Y\setminus K_h)+\|dT_h\|(Y\setminus K_h)<\epsilon\;.$$
Then, there exits a subsequence $(T_{h(n)})$ $w^*$-converging to $T\in N_k(Y)$. Moreover, $T$ has compact support if $\supp T_h$ are equi-bounded.\end{Teo}

\begin{Def}We say that $f:\R^k\to E$ is \emph{metrically differentiable} at $x\in\R^k$ if ther exists a seminorm $\|\cdot\|_x$ in $\R^k$ such that
$$\|f(x)-f(y)\|_E-\|x-y\|_x=o(|x-y|)\;.$$
This seminorm will be said to be the metric differential of $f$ at $x$ and denoted by $mdf(x)$.\end{Def}
\begin{Def}We say that $f:\R^k\to E$ is \emph{$w^*$-differentiable} at $x$ if there exists a linear map $L:\R^k\to E$ satisfying
$$w^*-\lim_{y\to x}\frac{f(x)-f(y)-L(y-x)}{|x-y|}=0\;.$$
This map $L$ will be said to be the $w^*$-differential of $f$ at $x$ and denoted by $wdf_x$.\end{Def}

\begin{Teo}For a Lipschitz map $f:\R^k\to E$, we have
$$mdf_x(v)=\|wdf_x(v)\|$$
for $\H^k-$a.e. $x\in \R^k$, for every $v\in\R^k$.\end{Teo}

\begin{Def}Given $S$ a countably $\H^k-$rectifiable set in $E$, if $f_i, A_i$ are as in Definition \ref{def_rect_set}, the \emph{approximate tangent space} to $S$ in $f_i(x)$ is
$$\mathrm{Tan}^{(k)}(S,f_i(x))=wd(f_i)_x(\R^k)\;.$$\end{Def}

This definition makes sense for $\H^k-$a.e. $x\in A_i$ and it is well posed as the dimension of the tangent is $\H^k-$a.e. $k$ and it does not depend on the $f_i'$s; moreover
$$\mathrm{Tan}^{(k)}(S_1,y)=\mathrm{Tan}^{(k)}(S_2, y)$$
for $\H^k-$a.e. $y\in S_1\cap S_2$.

\begin{Def}Let $L:\R^k\to E$ be linear. The \emph{$k$-jacobian} of $L$ is defined by
$$\mathbf{J}_k(L)=\frac{\omega_k}{\H^k(\{x\ :\ \|L(x)\|\leq1\})}=\frac{\H^k(\{L(x)\ :\ x\in B_1\})}{\omega_k}$$
where $B_1$ is the unit ball of $\R^k$ and $\omega_k=\H^k(B_1)$.\end{Def}
$\mathbf{J}_k$ satisfies the product rule for jacobians:
$$\mathbf{J}_k(L\circ M)=\mathbf{J}_k(L)\mathbf{J}_k(M)\;.$$
In the same way we can define the $k-$jacobian of a seminorm $s$ on $\R^k$:
$$\mathbf{J}_k(s)=\frac{\omega_k}{\H^k(\{x\ :\ s(x)\leq1\})}\;.$$

\begin{Teo}Let $f:\R^k\to E$ be a Lipschitz function. Then
$$\int_{\R^k}\theta(x)\mathbf{J}_k(mdf_x)dx=\int_E\sum_{x\in f^{-1}(y)}\theta(x)d\H^{k}(y)$$
for any Borel function $\theta:\R^k\to[0,+\infty]$ and
$$\int_A\theta(f(x))\mathbf{J}_k(mdf_x)dx=\int_E\theta(y)\H^0(A\cap f^{-1}(y))d\H^k(y)$$
for $A$ a Borel set in $\R^k$ and $\theta:E\to[0,+\infty]$ a Borel function.\end{Teo}

Moreover, we can also define the tangential differential on a rectifiable set.

\begin{Teo}Let $Y$ and $Z$ be duals of separable Banach spaces and let $S\subset Y$ be countably $\H^k-$rectifiable and let $g\in Lip(S,Z)$. Let $\theta:S\to(0,+\infty)$ be integrable with respect to $\H^k\llcorner S$ and set $\mu=\theta\H^k\llcorner S$.

Then, for $\H^k-$a.e. $x\in S$ there exis a linear and $w^*-$continuous map $L:Y\to Z$ and a Borel se $S^x\subset S$ such that $\Theta_k^*(\mu\llcorner S^x,x)=0$ (the upper $k-$dimensional density of $\mu\llcorner S^x$ in $x$) and
$$\lim_{S\setminus S^x\ni y\to x}\frac{d_w(g(y), g(x)+L(y-x))}{|y-x|}=0\;.$$
The map $L$ is uniquely determined on $\mathrm{Tan}^{(k)}(S,x)$ and its restriction to this space, denoted by $d^Sg_x$, satisfies the chain rule
$$wd(g\circ h)_y=d^Sg_{h(y)}\circ wdh_y\qquad \mathrm{for}\ \EL^k-\mathrm{a.e.}\ y\in A$$
for any Lipschitz function $h:A\to S$, $A\subset\R^k$.\end{Teo}

We have therefore a general area formula.

\begin{Teo}\label{teo_area}Let $g:E\to F$ be a Lipschitz function between Banach spaces and let $S\subset E$ be a countably $\H^k-$rectifiable set. Then
$$\int_S\theta(x)\mathbf{J}_k(d^Sg_x)d\H^k(x)=\int_F\sum_{x\in S\cap g^{-1}(y)}\theta(x)d\H^k(y)$$
for any Borel function $\theta:S\to[0,+\infty]$ and
$$\int_A\theta(g(x))\mathbf{J}_k(d^Sg_x)d\H^k(x)=\int_F\theta(y)\H^0(A\cap g^{-1}(y))d\H^k(y)$$
for any Borel set $A\subset E$ and any Borel function $\theta:F\to[0,+\infty]$.\end{Teo}

Finally, we mention a compactness theorem proved in \cite{ambrosio3}.

\begin{Teo}\label{teo_schmidt}Consider a normed space $E$ such that $E^*$ is separable, $n\in\N$, and a sequence $(T_h)\subset N_k(E^*)$ such that we have
$$M=\sup_hM(T_h)<+\infty\qquad M_d=\sup_hM(dT_h)<+\infty$$
and the $w^*-$tightness condition
$$\lim_{R\to\infty}\sup_{h}[\|T_h\|(E^*\setminus B_R(0))+\|dT_h\|(E^*\setminus B_R(0))]=0\;.$$
Then there exists a subsequence $(T_{h(n)})_n$ $w^*-$converging to $T\in N_k(E^*)$ with $M(T)\leq M$ and $M(dT)\leq M_d$.
\end{Teo}

\section{Metric currents on complex Banach spaces}

The results collected here and in the next subsection appear in \cite{mongodi1}, where we refer to for detailed proofs.

We examine the behavior of metric currents in relation with their projections on finite dimensional subspaces. In order to recover informations on the whole space from its finite dimensional subspaces, we introduce the following category of Banach spaces (see also \cite{noverraz1}).

A (complex) Banach space $E$  is said to have the \emph{projective approximation property (PAP for short)} if there exist a constant $a$ and an increasing collection $\{E_t\}_{t\in T}$ of finite dimensional subspaces of $E$ such that
\begin{enumerate}
\item[1)] $\{E_t\}_{t\in I}$ is a directed set for the inclusion;
\item [2)]$\displaystyle{E=\overline{\bigcup_{t\in I} E_t}}$;
\item[3)] for every $t\in I$ there exists a projection $p_t:E\to E_t$ with $\|p_t\|\leq a$.
\end{enumerate}

Every Banach space with a Schauder basis has the PAP. Two important cases of PAP Banach spaces with no Schauder basis are $\Ci(K)$, the space of continuous functions on a compact space K with the sup norm and $L^p(X,\mu)$, with $1\leq p\leq+\infty$, where $X$ is a locally compact space and $\mu$ being a positive Radon measure. In this section, we will work with Banach spaces having the PAP; we will endow the set $I$ of indeces with the partial ordering coming from the inclusion relation between subspaces.

\medskip

\begin{Prp}\label{prp_limlip}Let $f\in \Lip(E)$ and define $f_t=f\circ p_t$. Then $f_t\to f$ pointwise and $\Lip(f_t)\leq a\Lip(f)$, for every $t\in I$.
\end{Prp}

\begin{Prp}\label{prp_conv_proj}Let $T\in M_k(E)$ and define $T_t=(\pi_t)_\sharp(T)\in M_k(E_t)$ for every $t\in I$ such that $\dim_{\C}E_t\geq k$. By means of the inclusion $i_t:E_t\to E$, we can consider $T_t$ as an element of $M_k(E)$ and then, $T_t\to T$ weakly.\end{Prp}

\medskip

Let $\{E_t, p_t\}_{t\in I}$ be the countable collection of subspaces and projections given by PAP. We call it a \emph{projective approximating sequence} (PAS for short) if $p_t\circ p_s=p_{\min\{s,t\}}$.

We note that every separable Hilbert space or, more generally, every Banach space with a Schauder basis contains a PAS.

\begin{Teo}\label{teo_ext}Let us suppose that $\{E_t, p_t\}$ is a PAS in $E$. If we are given a collection of metric currents $\{T_t\}_{t\in I}$ such that 
\begin{enumerate}
\item $T_t\in N_k(E_t)$,
\item $(p_{t}\vert_{E_{t'}})_\sharp T_{t'}= T_{t}$ for every $t, t'\in I$ with $t'>t$,
\item $\|T_t\|\leq (p_t)_*\mu$ and $\|d T_t\|\leq (p_t)_*\nu$ for every $t\in I$ and some $\mu,\ \nu$ finite Radon measures on $E$. 
\end{enumerate}
then there exists $T\in N_k(E)$ such that $(p_t)_\sharp T=T_t$ for every $t\in I$. \end{Teo}

\medskip

We can substitute the request of the existence of a PAS and the compatibility condition (hypothesis \emph{(2)}) with an assumption on the existence of a global object. A \emph{metric functional} is a function $T:\E^k(E)\to\C$ which is subadditive and positively $1-$homogeneous with respect to every variable. For metric functionals, we can define mass, boundary and pushforward (see Section 2 of \cite{ambrosio1}).

\begin{Prp}Let $E$ be a Banach space with PAP. Suppose that $T:\E^k(E)\to\C$ is a metric functional with $T$ and ${\rm d}T$ of finite mass, such that $(p_t)_\sharp T\in N_k(E_t)$ for every $t\in I$. Then $T\in N_k(E)$.\end{Prp}

\subsection{Bidimension}

We say that $T\in\Di_m(U)$ is of bidimension $(p,q)$ if
$$T(f,\pi_1,\ldots,\pi_m)=0$$
whenever there exists $I\subset\{1,\ldots, m\}$, with $|I|>p$, such that $\pi_i\vert_{\supp(f)}$ is holomorphic for every $i\in I$, or $J\subset\{1,\ldots, m\}$, with $|J|>q$, such that $\pi_i\vert_{\supp(f)}$ is antiholomorphic for every $j\in J$. For a careful analysis of the notion of holomorphy in this context we refer the interested reader to the first chapters in \cite{noverraz1}. Here we only notice that Lipschitz holomorphic functions are not necessarily dense in the space of Lipschitz functions.

However, inspired by the links we found between the finite dimensional projections of a current and the current itself, we would like to give a different characterization of $(p,q)-$currents. 

We say that $T\in M_{k}(E)$ is \emph{finitely} of bidimension $(p,q)$ if every finite dimensional projection of it is a $(p,q)-$current.

\begin{Prp} $T\in M_k(E)$ is a $(p,q)-$current if and only if it is finitely so.\end{Prp}

As an application of Theorem \ref{teo_ext}, we have the following result about the existence of a Dolbeault decomposition for $T\in M_k(E)$.

\begin{Prp}\label{prp_dec}Let us suppose that $\{E_t, p_t\}$ is a PAS in $E$. Let $T\in N_k(E)$; if $T_t$ has a Dolbeault decomposition in normal $(p,q)-$currents in $E_t$ for all $t\in T$, with a finite Radon measure $\nu$ (independent of $t$) whose pushforward dominates the boundaries' masses, then also $T$ admits a Dolbeault decomposition.\end{Prp}

\begin{Rem} In general it is not easy to verify the hypotheses of Proposition \ref{prp_dec} for a current $T\in N_k(E)$; however, this result is an example of a general phenomenon: in a Banach space with the projective approximation property, it is often enough to check a certain property for finite dimensional subspaces in order to obtain that it holds for the whole space. For instance, any equality between currents holds in $E$ if and only if it holds finitely, namely whenever the currents are pushed forward through a finite rank projection.\end{Rem}

Employing the idea given in this Remark, we can show the following.

\begin{Cor}If $T\in M_k(E)$ admits a Dolbeault decomposition, then it is unique.\end{Cor}

Let us suppose that $T$ is a $(p,q)-$current whose boundary admits a Dolbeault decomposition. Then we can define $\de T$ and $\debar T$ as follows.

Write $dT=S_1+\ldots+S_h$ with $S_i\in M_{p_i,q_i}(U)$ where $p_i+q_i=p+q-1=m$ since $dT\in M_{p+q-1}(U)$. If $(f,\pi)$ is a $m-$form of pure type $(p_i,q_i)$, then
$$S_{p_i,q_i}(f,\pi)=dT(f,\pi)=T(1,f,\pi)$$
if $p_i>p$ or $q_i>q$,and consequently $T(1,f,\pi)=0$, since $T$ is a $(p,q)-$current. Therefore, we can only have two cases: $p=p_i$ and $q-1=q_i$ or $p-1=p_i$ and $q=q_i$ i.e.
$$dT=S_{p,q-1}+S_{p-1,q}$$
and we put
$$\de T=S_{p-1,q}\qquad \debar T=S_{p,q-1}$$
Therefore, if a current $T$ admits a decomposition in $(p,q)$ components, we can define $\de T$ and $\debar T$ setting
$$
\de T=\sum_{i=1}^h\de T_i\qquad \debar T=\sum_{i=1}^h\debar T_i
$$
where $T=T_1+\ldots+T_h$ is the $(p,q)$ decomposition.

\medskip

\begin{Prp}\label{prp_pushf_debar}
If $H:E\to F$ is a holomorphic map between complex Banach spaces, then, for every current $T\in \D_m(E)$ for which $\de T$ and $\debar T$ are defined, the following hold:
$$H_\sharp \de T=\de H_\sharp T\qquad H_\sharp\debar T=\debar H_\sharp T$$
\end{Prp}

Moreover, it is easy to check that $C_m\circ\de=\de\circ C_m$ and $C_m\circ\debar=\debar\circ C_m$, where $C_m$ is the map given by Theorem \ref{teo_comp_man}.

\begin{Prp}\label{prp_formule_de}
We have that $\de^2=\debar^2=0$ and $\de\debar=-\debar\de$
\end{Prp}

\subsection{$\debar$ and slicing}

Let $T$ be a normal metric current whose boundary admits a Dolbeault decomposition. Then $\supp \debar T, \supp \de T\subseteq \supp dT$ and, moreover, $\|\debar T\|_A, \|\de T\|_A\leq C\| dT\|_A$ for every $A\subset X$. In particular, if $dT$ is rectifiable, then $\debar T$ and $\de T$ are too.

\medskip

The slices of a current $T$ through a map $\pi:E\to\R^{n}$ are defined by
$$\langle T,\pi, x\rangle(f,\eta)=\lim_{\epsilon\to 0} T(\rho_{\epsilon, x}f, \pi, \eta)=\lim_{\epsilon\to 0}(-1)^{n(k-1-n)}\pi_\sharp(T\llcorner(f,\eta))(\rho_{\epsilon, x}, x_1,\ldots, x_n)$$
with $\rho_{\epsilon, x}$ any family of smooth approximations of  $\delta_x$. If $\pi:E\to\R^{2n}\cong\C^n$, we can write the slices as
$$\langle T,\pi, x\rangle=\lim_{\epsilon\to 0} \pi_\sharp(T\llcorner(f,\eta))(\rho_{\epsilon, x}, z_1,\bar{z}_1,\ldots, z_n,\bar{z}_n)\;.$$

\begin{Prp}The operators $\de$ and $\debar$ commute with the slicing through holomorphic maps.\end{Prp}
\noindent{\bf Proof: }Let $T\in N_{k}(E)$ be a  normal current, such that $\de T$ and $\debar T$ are again  normal. Let $\pi:E\to\C^n$ be a holomorphic map.

\medskip

We start with the finite dimensional case, supposing that $E=\C^N$.

Let $(f,\eta)$ be a $(k-2n-1)-$metric form with $\Ci^2$ coefficients. We treat only the case of $\debar T$, the proof for $\de T$ being analogous.

Let $z_1,\ldots, z_n$ and $w_1,\ldots, w_N$ be holomorphic coordinates in $\C^n$ and $\C^N$, respectively. The slice $\langle \debar T, \pi, x\rangle$ exists for a.e. $x$, by Theorem \ref{teo_ext_slice} and we have
$$\langle \debar T,\pi, x\rangle(f,\eta)=\lim_{\epsilon\to 0}\pi_\sharp((\debar T)\llcorner(f,\eta))(\rho_{\epsilon,x}, z_1,\bar{z}_1,\ldots, z_n, \bar{z}_n)\;.$$
Now, we set  $\widetilde{\eta}^{jh}$ to be the $(k-2n-1)-$tuple differing from $\eta$ only in the $h-$th component, which is $(-1)^h\de \eta_h/\de\bar{w}_j$. By \cite{mongodi1}[Proposition 6]
$$(\debar T)\llcorner(f,\eta)=(-1)^{k-1}\debar(T\llcorner (f,\eta)) + \sum_{j=1}^NT\llcorner(\de f/\de \bar{w}_j, \bar{w}_j, \eta)+\sum_{h,j}T\llcorner(f, \bar{w}_j,\widetilde{\eta}^{jh})$$
and we note that $T\llcorner(f,\eta)$ is a $2n+1-$form, so 
$$\pi_\sharp\debar T\llcorner(f,\eta)=\debar\pi_\sharp T\llcorner(f,\eta)=0$$
by Proposition \ref{prp_pushf_debar}, as $\pi$ is holomorphic. It follows
$$\pi_\sharp((\debar T)\llcorner(f,\eta))(\rho_{\epsilon, x},z_1,\ldots, \bar{z}_n)=\sum_{j=1}^N\pi_\sharp(T\llcorner(\de f/\de\bar{w}_j, \bar{w}_j, \eta))(\rho_{\epsilon, x},z_1,\ldots, \bar{z}_n)+$$
$$+\sum_{j,h}\pi_\sharp(T\llcorner(f,\bar{w}_j,\widetilde{\eta}^{jh}))(\rho_{\epsilon,x},z_1,\ldots,\bar{ z}_n)$$
$$=\sum_{j=1}^NT\left(\frac{\de f}{\de \bar{w}_j}\cdot(\rho_{\epsilon,x}\circ\pi), \pi_1, \bar{\pi}_1,\ldots, \pi_n,\bar{\pi}_n, \eta\right)+\sum_{j,h}T\left(f\cdot(\rho_{\epsilon,x}\circ\pi), \pi_1,\ldots, \bar{\pi}_n, \widetilde{\eta}^{jh}\right)$$
$$=\sum_{j=1}^NT\llcorner(\rho_{\epsilon,x}\circ\pi, \pi_1,\ldots, \bar{\pi}_n)(\de f/\de\bar{w}_j,\bar{w_j},\eta)+\sum_{j,h}T\llcorner(\rho_{\epsilon,x}\circ\pi,\pi_1,\ldots,\bar{\pi}_n)(f,\bar{w}_j,\widetilde{\eta}^{jh})\;.$$
Now, again by \cite{mongodi1}[Lemma 1], we have
$$\debar\langle T,\pi, x\rangle(f,\eta)=\sum_{j=1}^N\langle T,\pi, x\rangle\left(\frac{\de f}{\de\bar{w}_j}, \bar{w}_j, \eta\right)+\sum_{j,h}\langle T,\pi,x\rangle(f,\bar{w}_j,\widetilde{\eta}^{jh})$$
$$=\sum_{j=1}^N\lim_{\epsilon\to0} T\llcorner(\rho_{\epsilon,x}\circ\pi, \pi_1,\ldots, \bar{\pi}_n)\left(\frac{\de f}{\de\bar{w}_j}, \bar{w}_j,\eta\right)+\sum_{j,h}\lim_{\epsilon\to0} T\llcorner(\rho_{\epsilon,x}\circ\pi, \pi_1,\ldots, \bar{\pi}_n)\left(f, \bar{w}_j,\widetilde{\eta}^{jh}\right)$$
and by the previous computation this is equal to
$$\lim_{\epsilon\to 0}\pi_{\sharp}((\debar T)\llcorner(f,\eta))(\rho_{\epsilon, x}, \pi_1,\ldots,\bar{\pi}_n)=\langle \debar T,\pi, x\rangle(f,\eta)\;.$$
Therefore, for every $(f,\eta)$ with $f\in\Ci^2$
$$\langle \debar T,\pi, x\rangle(f,\eta)=\debar\langle T,\pi, x\rangle (f,\eta)\;.$$
This means that the functional $\debar\langle T,\pi, x\rangle$ can be extended to all the metric forms as a metric current, by defining it equal to $\langle \debar T,\pi, x\rangle$. 

\medskip

If $E$ is infinite-dimensional, the previous argument gives that 
$$\langle \debar T,\pi, x\rangle(f,\eta)=\debar\langle T,\pi, x\rangle (f,\eta)$$
for every metric form whose coefficients depend on a finite number of variables. This implies that
$$(p_t)_\sharp\langle \debar T,\pi, x\rangle=(p_t)_\sharp\debar\langle T,\pi, x\rangle (f,\eta)$$
for every finite rank projection $p_t$, and consequently that 
$$\langle \debar T,\pi, x\rangle=\debar\langle T,\pi, x\rangle\;. $$
$\Box$

\medskip

As already said, the previous proof works also for the $\de$:
$$\de\langle T,\pi, x\rangle=\langle \de T,\pi, x\rangle\;.$$
Moreover, if $\de\debar T$ is a metric current, then
$$\de\debar\langle T,\pi, x\rangle=\langle \de\debar T,\pi, x\rangle\;,$$
or, which is the same,
$$dd^c\langle T,\pi,x\rangle=\langle dd^c T,\pi, x\rangle\;.$$

\section{Analytic sets in Hilbert spaces}

There are many possible definitions for a finite dimensional analytic set in an infinite-dimensional space. Here we adopt the following (see \cite{ruget1}).
 \begin{Def}Let $H$ be a complex Hilbert space (or more generally a Banach space). A closed set $A\subset H$ will be called a {\em finite-dimensional analytic set} in $H$ if, locally in $H$, $A$ is an analytic subspace of a complex submanifold of $H$ of finite dimension .\end{Def}

\begin{Rem} This definition is equivalent to the one given by Douady in \cite{douady1} (see \cite{ruget1} for the details).\end{Rem}

An equivalent definition is the following (see \cite{aurich1}): $A\subset H$ is a finite-dimensional analytic set if for each $v\in H$ there exist a neighborhood $U$, another complex Hilbert space $H'$ and a holomorphic map $F:U\to H'$ whose differential has finite dimensional kernel such that $A\cap U=F^{-1}(0)$. This follows easily from the Implicit Function Theorem; this turns out to be equivalent to asking for a map whose differential has finite dimensional kernel and cokernel (i.e. a Fredholm map) such that $A\cap U=F^{-1}(0)$. 

\medskip

We recall a result from \cite{ruget1}.

\begin{Teo}\label{teo_hilb_img} Let $X$ be an analytic subset of finite dimension, $W$ be a hilbertian manifold and $f:X\to W$ a proper holomorphic map. Then, $f(X)$ is a finite-dimensional analytic subset of $W$.\end{Teo}

\noindent{\bf Example} The properness assumption cannot be dropped as shown by the following example. Let us consider the space $\ell^2$ of square-summable sequences of complex numbers and consider the holomorphic map $g:\mathbb{D}^2\to\ell^2$ given by
$$g(z,w)=\{zw^n\}_{n\in\N}\;.$$
The preimage of $\{0\}$ is $\{z=0\}$, which is not compact. Let $X=g(\DD^2)$ and assume, by a contradiction, that there exists a holomorphic function $\Phi:\ell^2\to\C$ is a holomorphic function vanishing on $X$; then $\Phi\circ g=0$ and consequently. 
$$0=\Phi\circ g(z,w)=\sum_{n,m\geq0}\frac{\de^{n+m} \Phi\circ g}{\de z^n\de w^m}(0,0)\frac{z^nw^m}{n!m!}$$
One has
$$0=\frac{\de \Phi\circ g}{\de z}(0,0)=\frac{\de\Phi}{\de e_0}(0)$$
$$0=\frac{\de^2\Phi\circ g}{\de w\de z}(0,0)=\frac{\de}{\de w}\left(\sum_j\frac{\de \Phi}{\de e_j}(g)\frac{\de g_j}{\de z}\right)(0,0)=$$
$$\left(\sum_{j,k}\frac{\de^2 \Phi}{\de e_k\de e_j}(g)\frac{\de g_k}{\de w}\frac{\de g_j}{\de z} + \sum_j\frac{\de\Phi}{\de e_j}\frac{\de^2 g_{j}}{\de w\de z}\right)(0,0)=\frac{\de\Phi}{\de e_1}(0,0)\;,$$
since
$$\frac{\de g_k}{\de w}(0,0)=0\qquad \forall k$$
and
$$\frac{\de^2g_k}{\de z\de w}(0,0)\neq 0 \Leftrightarrow g_k(z,w)=zw \Leftrightarrow k=1\;.$$
Proceeding in this way, we can show that $\de\Phi/\de e_j=0$ in $0\in \ell^2$, for all $j\in\N$. Therefore all the derivatives of $\Phi$ vanish at the origin, which means that no regular hypersurface of $\ell^2$ can contain a neighborhood of $0\in X$.

\bigskip

The major advantage of the given definition is that we can recover the local properties of finite-dimensional analytic sets in an infinite dimensional space from the usual local results on analytic sets in a complex manifold. In particular, let $X$ be a finite-dimensional analytic set in $H$, then
\begin{enumerate}
\item $X$ admits a local decomposition in finitely many irreducible components;
\item such a decomposition is given by the closure of the connected components of the regular part of $X$;
\item for every $x\in X$ there exist a finite-dimensional subspace $L$ of $H$ and an orthogonal projection $\pi:H\to L$ which realizes a neighborhood of $x\in X$ as a finite covering on $L$;
\item $X$ is locally connected by analytic discs;
\item if $X$ is irreducible, every nonconstant holomorphic function is open;
\item if $X$ is irreducible, the maximum principle holds;
\item if $X$ is compact and $H$ is holomorphically separable, then $X$ is finite.
\end{enumerate}

\bigskip

The behaviour of the analytic sets in a Banach space can vary wildely, depending on the properties of the space. We give here some examples.

\medskip

\noindent{\bf Example } Let $c_0$ be the vector space of sequences of complex numbers vanishing at infinity, i.e. $\{a_n\}\subset\C$ such that $\lim_{n\to\infty}a_n=0$; $c_0$ is a Banach space with the supremum norm. We consider the holomorphic map $f:\DD\to c_0$ given by
$$f(z)=\{z^n\}_{n\in\N}\;;$$ $f$ is a regular and injective holomorphic map; its image is contained in the unit ball of $c_0$ and if $\{z_k\}$ is a sequence converging to $b\DD$,  $\{f(z_k)\}_k$ does not converge, therefore $f$ is proper. Thus, $f(\DD)$ is a complex manifold of dimension $1$ in $E$, which is bounded.

\medskip

\noindent{\bf Example } Generalizing the previous example, we consider the Banach space of $p-$summable sequences of complex numbers $\ell^p$ and the holomorphic map $F:\DD^k\to\ell^p$ given by
$$F(z_1,\ldots, z_k)=\{z^I\}_{I}$$
where $I$ varies through all the multi-indeces of length $k$. We have that
$$|z^I|\leq (\max\{|z_1|,\ldots, |z_k|\})^{|I|}$$
and the number of multiindexes $I$ with $|I|=i_1+\ldots+i_k=m$ is
$${m+k-1\choose k-1}$$
which is less than $(2m)^k$ if $m$ is large enough. Therefore
$$\sum_I |z^I|^p\leq\sum_{m} (2m)^k(\max\{|z_1|,\ldots, |z_k|\})^m$$
which converges for $\max\{|z_1|,\ldots, |z_k|\}<1$.

Again, the map $F$ is regular, injective and proper with unbounded image $F(\DD^k)$. We observe that $F(\DD^k)$ provides an example of finite dimensional manifold not contained in any finite-dimensional linear subspace.

\subsection{Positive currents}

\begin{Def}Let $T$ be a metric $(p,p)-$current. We say that $T$ is \emph{positive} if, given $\pi_1,\ldots, \pi_p\in\Ol(X)$, 
$$T(f,\pi_1,\overline{\pi}_1,\ldots, \pi_p,\overline{\pi}_p)\geq0$$ for every Lipschitz function $f\geq0$ on $X$. We say that a current is \emph{finitely positive} if every finite dimensional projection of it is positive.\end{Def}

\begin{Prp}Let $\Omega\subset H$ be an open set. $T\in M_{2p}(\Omega)$ is positive if and only if it is finitely positive.\end{Prp}
\noindent{\bf Proof: } Obviously, if $T$ is positive then every complex linear pushforward of it is positive. On the other hand, if $p_m:H\to\C^m$ is the projection on the first $m$ coordinates, then, by Proposition \ref{prp_limlip}, $f\circ p_m\to f$ pointwise and $\Lip(f\circ p_m)\leq\Lip(f)$.

Therefore, 
$$T(f,\pi_1,\ldots,\overline{\pi}_p)=\lim_{m\to\infty}T(f\circ p_m, \pi_1\circ p_m,\ldots, \overline{\pi}_p\circ p_m)=$$
$$\lim_{m\to\infty}(p_m)_\sharp T(f\vert_{p_m(H)},\pi_1\vert_{p_m(H)},\ldots, \overline{\pi}_p\vert_{p_m(H)})\geq0\;,$$
which is the thesis. $\Box$

\medskip

On an infinite dimensional complex space,the plurisubharmonic functions are defined by the submean property: let $\Omega\subset H$ be an open set and let $u:\Omega\to\R\cup\{-\infty\}$ be an upper semicontinuous function (not identically equal to $-\infty$); $u$ is plurisubharmonic if
$$u(a)\leq\int_{0}^{2\pi}u(a+e^{i\theta}b)d\theta$$
for every $a\in\Omega$ and every $b\in H$ such that $a+\lambda b\in\Omega$ for every $\lambda\in\C$ with $|\lambda|\leq1$. See also \cite{mujica1} for a discussion of the properties of such functions.

Define $d^c=i(\de-\debar)$.

\begin{Prp}Let $T\in M_{2p}(\Omega)$ be a positive closed current with bounded support and $u:\Omega\to\R\cup\{-\infty\}$ a bounded plurisubharmonic function. Then $dd^c(T\llcorner u)$ is a closed, positive metric current with bounded support and the following estimate holds:
$$M(dd^c(T\llcorner u))\leq \|u\|_\infty M(T)\;.$$
\end{Prp}
\noindent{\bf Proof: } We note that the result is true for any finite-dimensional projection of $T$. Namely, if $p_m$ is as above, 
$$(p_m)_\sharp dd^c(T\llcorner u)=dd^c((p_m)_\sharp T\llcorner (u\circ p_m))=dd^c(T_m\llcorner u_m)$$
with $T_m$ a positive closed current with compact support in $\C^m$ and $u_m$ a bounded plurisubharmonic function. Then we know that $$dd^c(T_m\llcorner u_m)=T_m\llcorner (dd^c u_m)$$
in the sense of distributions and for every compact $K$  we have
$$M_K(dd^c(T_m\llcorner u_m))\leq C \|u_m\|_{\infty, K}M_K(T_m)\;.$$
As $p_m$ is an orthogonal projection, $M_B(T_m)\leq M_B(T)$, whereas the fact that $T_m$ converges weakly to $T$, together with the semicontinuity of the mass, implies that $M_B(T)\leq\limsup M_B(T_m)$, therefore $M_B(T_m)\to M_B(T)$ for every bounded set $B$ in $\Omega$.

This means that, for $j$ big enough,
$$M_K(dd^c(T_{m_j}\llcorner u_{m_j}))\leq C\|u\|_{\infty, K}M_K(T)$$
so,  by \ref{teo_schmidt}, we can find a subsequence $w*-$converging to some $S\in M_{2p-2}(\Omega)$.

Such an $S$ is such that $(\pi_m)_\sharp S=T_m$ for infinitely many $m$, therefore $S$ coincides, as a metric functional, with $dd^c(T\llcorner u)$. Moreover, $S$ is positive and closed and
$$M_K(S)\leq \|u\|_{\infty,K}M_K(T)\;.$$
We note also that $\supp S\subseteq \supp T$. $\Box$

\medskip

In view of the previous Proposition, we will denote by $dd^cu\wedge T$ the current $dd^c(T\llcorner u)$.

Proceeding by induction, we can give a meaning to the writing
$$dd^cu_1\wedge\ldots\wedge dd^c u_p\wedge T$$
for a current $T$ in the hypotheses of the previous Proposition. Such a definition allows us to write an analogue of the Monge-Amp\`ere operator in Hilbert spaces.

\subsection{Currents of integration on analytic sets}

Let $V$ be a finite-dimensional analytic set in some open domain $U\subset H$; since, by definition, $V$ is locally contained in some finite-dimensional submanifold, we know that it is locally of finite volume. Therefore, for any $p\in V$ there exists a ball $B$ such that the current $[V]\llcorner B$ of integration on the regular part of $V\cap B$ is a well-defined rectifiable metric current. The following result gives an estimate for the mass of such a current, analogous to Wirtinger formula in the finite dimensional case.

\begin{Prp}\label{prp_vol_in_H}Let $H$ be a Hilbert space, with scalar product $\langle\cdot,\cdot\rangle$, $V$  an analytic set in an open set $U\subset H$, with $\dim_\C V_\rg=p$. Let $\Omega$ be a ball in $U$ and let $[V]$ be the current of integration associated to $V\cap \Omega$ in $\Omega$. Then
$$M([V])\leq\lim_{n\to\infty}\sum_{1\leq i_1<\ldots<i_p\leq n}[V]\left(\frac{i^p}{2^pp!}, z_{i_1},\bar{z}_{i_1},\ldots, z_{i_p},\bar{z}_{i_p}\right)<+\infty\;,$$
where $\{z_j\}_{j\in \N}$ are coordinate functions with respect to some orthonormal basis.
\end{Prp}
\noindent{\bf Proof: }
Given an orthonormal basis $\{e_n\}_{n\in\N}$, let $E_m=\mathrm{span}\{e_1,\ldots, e_m\}$ and $\pi_m:H\to E_m$ the orthogonal projection. We denote $[V]_m=(\pi_m)_\sharp [V]$ and observe that $[V]_m\to [V]$ weakly, so that by the semicontinuity of the mass we get
$$M([V])\leq \liminf_{m\to\infty}M([V]_m)\;,$$
where the masses are relative to $\Omega$ or to $\Omega_m=\Omega\cap E_m$ respectively.

On the other hand, the projections $\pi_m$ have norm $\|\pi_m\|\leq 1$, so, if $\mu$ is the mass measure of $[V]$, we have
$$|[V]_m(f,\eta)|\leq \prod \Lip(\eta_j)\int_{U_m}|f|d(\pi_m)_\sharp \mu\;.$$
This implies that the mass measure $\mu_m$ of $E_m$ is dominated by $(\pi_m)_\sharp \mu$, therefore
$$\mu_m\leq (\pi_m)_\sharp \mu\xrightarrow[m\to\infty]{}\mu$$
which means that 
$$M([V])=\mu(\Omega)\leq\liminf_{m\to\infty}\mu_m(\Omega)\leq \lim_{m\to\infty}(\pi_m)_\sharp \mu(\Omega)=\mu(\Omega)=M([V])\;.$$
Now, $E_m$ with the induced scalar product is the usual complex hermitian space $\C^m$ and the pushforward of an analytic chain is again an analytic chain. Therefore
$$M([V]_m)\leq[V]_m(\omega_m^p/p!)=\H^{2p}(V_m)\leq C_p M([V]_m)$$
where
$$\omega_m=\frac{i}{2}\sum_{j=1}^mdz_j\wedge d\bar{z}_j=\frac{i}{2}\sum_{j=1}^me^*_j\wedge \bar{e}^*_j$$
and $C_p$ is a constant depending only on $p$ (see \cite{ambrosio1}, after Remark 8.4).
Noticing that $[V]_m(\omega_m^p/p!)=[V](\omega_m^p/p!)$, we obtain
$$M([V])\leq \lim_{n\to\infty}\sum_{1\leq i_1<\ldots<i_p\leq n}[V]\left(\frac{i^p}{2^pp!}, z_{i_1},\bar{z}_{i_1},\ldots, z_{i_p},\bar{z}_{i_p}\right)\leq C_pM([V])<+\infty\;,$$
which is the thesis. $\Box$

\medskip

\begin{Rem}The current $[V]$ is obviously positive, of bidimension $(k,k)$ for some $k$ and its boundary is supported outside $\Omega$.\end{Rem}

Given an orthonormal basis $\{e_j\}_{j\in\N}$ and a multiindex $I=(i_1,\ldots, i_k)$, let $\pi_I$ denote the orthogonal projection from $H$ onto $\mathrm{Span}\{e_{i_1},\ldots, e_{i_k}\}$.

\begin{Teo}\label{teo_king}Let $\Omega\subset H$ be a ball, $S$ be a rectifiable current in $\Omega$. Suppose that
\begin{enumerate}
\item $\supp dS\cap\Omega=\emptyset$;
\item $S$ is a $(k,k)$ positive current.
\end{enumerate}
Then $S$ can be represented as a sum with integer coefficients of integrations on the regular parts of analytic sets. \end{Teo}

\begin{Rem}We cannot show that $X$ is a finite-dimensional analytic space in the sense precised in the beginning of this section; indeed, in the example discussed before, the map $f(z,w)=(zw^n)_{n}$, gives an analytic space which carries a current of integration which satisfies the hypotheses of the previous theorem but cannot be written ad the integration on a finite-dimensional analytic space.
\end{Rem}

\medskip

\noindent{\bf Proof: } Since $S$ is a metric current, we can define its pushforward through any Lipschitz map.  We note that $(\pi_I)_\sharp S=m_I[V_I]$ with $V_I=\pi_I(\Omega)$. By \emph{(ii)} we know that $m\geq0$ and by the fact that $dS\llcorner\Omega=0$, we deduce that $S$ is integral in $\Omega$, therefore $m\in\N$.

By Theorem \ref{teo_ext_slice}, for almost every $y\in V_I$ we can define $\langle S,\pi_I,y\rangle$; moreover, we can find a $\H^{2k}-$rectifiable subset $B$ of $\supp S$ and an integer multiplicity function $\theta(x)$ such that $S=[B]\llcorner\theta$; then
$$\langle S,\pi_I,y\rangle=\sum_{x\in\pi_I^{-1}(y)\cap B}\theta(x)[x]$$
and
$$\sum_{x\in\pi_I^{-1}(y)\cap B}\theta(x)=\langle S,\pi_I, y\rangle(1)=m\;.$$
Let us call $G_I\subset V_I$ the set of $y$ such that the slice exists; then for $j\not\in I$ and $z\in G$, we set
$$P_j(z,W)=\prod_{x\in\pi_I^{-1}(z)\cap B}(W-w_j(x))^{\theta(x)}$$
where $w_j(x)$ is the $j-$th coordinate of $x$ in the fixed orthonormal basis.

We note that
$$\sum_{x\in\pi_I^{-1}(z)\cap B}\theta(x)w_j(x)^s=\langle S,\pi_I, z\rangle(w_j^s)$$
is a holomorphic function of $z$, because $\debar S=0$, therefore by a classical argument $P_j(z,W)$ is a polynomial with coefficients in $\Ol(V_I)$ for every $j\not\in I$.

After removing an $\H^{2k}-$negligible set from $G_I$, we have that $P_j(z,w_j)=0$ for every $j\not\in I$ and every $x=(z,w)\in \pi_I^{-1}(G_I)\cap B$.

Let us define
$$X_I=\{P_j(z,w_j)=0\ ,\ j\not\in I\}\qquad X=\bigcup_I X_I\;.$$
We can look at $X_I$ as the zero locus of the map 
$$P_I:H\to \mathrm{Span}_\C\{e_j\ :\ j\not\in I\}=H_1$$
given by
$$P_I(z,w)=\sum_{j\not\in I}e_jP_j(z,w_j)\;.$$
In order to show that $P_I$ is well defined, we observe that, since $P_j(z,W)$ is a polynomial in $W$, for a fixed $z$ we have 
$$|P_j(z,W)|\leq \min\{d(W, w_j(x))^m,\ x\in\pi_I^{-1}(z)\cap B\}$$
with $m\geq 1$, if $W$ is close enough to some $w_j(x)$. Therefore, for $p=(z,w)$ in a neighborhood of $B$, we can write
$$\sum_{j\not\in I}|P_j(z,w_j)|^2\leq\sum_{j}|w_j-w_j(x)|^{2m}\leq\|p-x\|^{2m}<+\infty$$
where $x$ is the nearest point in $B$ to $p=(z,w)$. The map $P_I$ is obviously holomorphic, as its entries are polynomials in $w_j$ with coefficients holomorphic in the first $k$ coordinates.

This shows that $X_I$ is locally given as the zero locus of a holomorphic map between Hilbert space, therefore it is an analytic set. Moreover, let $x\in X_I$ be a smooth point and suppose that $\dim_\C T_x X_I>k$. By construction, $\pi_I\vert_{T_xX_I}:T_xX_I\to\C^k$ has maximum rank, therefore we can find $m$ such that, setting $J=I\cup\{m\}$, the projection $\pi_J$ restricted to $T_xX_I$ is surjective onto $\C^{k+1}$. This means that $P_m(z,W)\equiv 0$, but this is impossible. So, $\dim_\C T_xX_I=k$, i.e. the regular part of $X_I$ is a smooth $k-$dimensional complex manifold.

Now, let
$$B_I=\{x\in B\ :\ \mathbf{J}_{2k}\pi_I(x)\neq 0\}$$
i.e. the set of points $x$ of $B$ such that $D\pi_I$ has rank $2k$ on the approximate tangent to $B$ at $x$; define also
$$C_I=(B\cap V) \setminus \pi_I^{-1}(G_I)\;.$$
By Theorem \ref{teo_area}, we have
$$\int_{B_I\cap C_I} \mathbf{J}_{2k}\pi_I(x)d\H^{2k}(x)=\int_{V_I\setminus G_I}\left(\int_{\pi_I^{-1}(y)\cap B} gdH^0\right)d\H^{2k}(y)=0$$
where $g$ is the characteristic function of $B_I\cap C_I$. Since $J_{2k}\pi_I>0$ on $B_I\cap C_I$, this means that $H^{2k}(B_I\cap C_I)=0$.

Obviously, $B=\bigcup B_I$ and
$$B\cap\left(\bigcup_I \pi_I^{-1}(G_I)\right)\subset X\;,$$
but
$$B\setminus \bigcup_I \pi_I^{-1}(G_I)=\bigcap_{I}(B\setminus \pi_I^{-1}(G_I))\subset \bigcap_{I}((B\setminus B_I)\cup(B_I\cap C_I))\subseteq\bigcup_I(B_I\cap C_I)=D\;.$$
Since $\H^{2k}(B_I\cap C_I)=0$, we also have $\H^{2k}(D)=0$ and $\H^{2k}(B\setminus X)=0$. Moreover, as $X$ is closed in $\Omega$, $\supp S\subset X$; therefore $S=[B\cap X]$.

\medskip

If we denote by $X_\rg$ the union of the regular parts of $X_I$, then $S\llcorner X_\rg$ is a $(k,k)-$current, positive and closed,  with support on a $k-$dimensional smooth complex manifold. Therefore, $S\llcorner X_\rg$ can be written as a series with integer coefficients of the currents of integration on the connected components of $X_\rg$.

There exists $r>0$ such that $\pi_I(V)$ contains a ball of radius $r$ for every $I$; therefore, the $\H^{2k}-$measure of the regular part of $X_I$ is uniformly bounded from below independetly of $I$. On the other hand $S\llcorner X_\rg$ is of finite mass; therefore it has to be a finite sum.

Finally, let us consider the rectifiable set $R=B\setminus X_\rg$. If we project it on the first $m$ coordinates, for $m\geq k+1$, we obtain that its image is the singular set of a $k-$dimensional analytic space, therefore $\H^{2k}-$negligible; again by  Theorem \ref{teo_area},
$$\int_{R}\mathbf{J}_{2k}({d^R\pi_m}_x)d\H^{2k}(x)=\int_{\C^m}\sharp\{x\in R\cap \pi_m^{-1}(y)\}d\H^{2k}$$
with $\pi_{m}:H\to\mathrm{Span}\{e_1,\ldots, e_m\}$ the orthogonal projection. Let us denote by $\eta(x)$ the approximate tangent to $R$ in $x$; then the $2k-$jacobian of $\pi_m$ on $\mathrm{Tan}^{(2k)}(R,x)$ is given by the projection of $\eta(x)$ on $\mathrm{Span}\{e_1,\ldots, e_m\}$.

We define $A_m=\{ x\in R\ :\ \mathbf{J}_{2k}({d^R\pi_m}_x)>0\}$ and we note that $A_{k+1}\cup A_{k+2}\cup\ldots=R$, up to $\H^{2k}-$neglibigle sets. But
$$\int\limits_{A_m}\mathbf{J}_{2k}({d^R\pi_m}_x)d\H^{2k}(x)=\int\limits_{\C^m}\sharp\{x\in A_m\cap \pi_m^{-1}(y)\}d\H^{2k}=\int\limits_{\pi_m(A_m)}\sharp\{x\in R\cap \pi_m^{-1}(y)\}d\H^{2k}=0$$
because $\pi_m(A_m)\subseteq \pi_m(R)$, which is $\H^{2k}-$negligible. Therefore $\H^{2k}(R)=0$, so $S\llcorner X_\rg=S$ and this concludes the proof. $\Box$

\medskip

\section{Boundaries of holomorphic chains}

\begin{Prp}\label{prp_emb}Let $M$ be a compact $\Ci^1$ submanifold of a complex reflexive Banach space $E$ with complex structure $J$, such that $\dim_\R M=2p-1$ and $\dim_\C T_zM\cap JT_zM=p-1$ for every $z\in M$. Then there exists a complex linear map $F:E\to\C^n$ for some $n>0$ which restricts to an embedding of $M$ into $\C^n$.\end{Prp}
\noindent{\bf Proof: }Given $z\in M$, let $l_{1,z},\ \ldots, l_{p,z}$ linearly independent elements of $E^*$ such that
$$\ker l_{1,z}\cap\ldots\cap\ker l_{p-1,z}\cap T_{z}M\cap JT_zM=\{0\}$$
and
$$\ker l_{1,z}\cap\ldots\cap\ker l_{p-1,z}\cap \ker\mathsf{Re}\ l_{p,z}\cap T_{z}M=\{0\}\;;$$
both these conditions are open. By compactness, we can find finitely many $l_1,\ldots, l_N$ such that
for every point $z\in M$ there exists indexes $j_1<\ldots<j_p$ such that
$$\ker l_{j_1}\cap\ldots\cap\ker l_{j_{p-1}}\cap T_{z}M\cap JT_zM=\{0\}$$
and (viewing $E$ as a real vector space)
$$\ker l_{j_1}\cap\ldots\cap\ker l_{j_{p-1}}\cap \ker\mathsf{Re}\ l_{j_p}\cap T_{z}M=\{0\}\;.$$
This means that if we define $L:E\to\C^N$ by $L=(l_1,\ldots, l_N)$, we have that the differential $dL$ is always of real rank $2p-1$ on $M$ and it is complex linear on $T_zM\cap JT_zM$.

\medskip

Let $U_1,\ldots, U_h$ be the open sets and $l_1,\ldots, l_p, l_{p+1},\ldots, l_{2p},\ldots, l_{hp}$ be the maps constructed as above and $\{V_j\}_{j=1}^K$ a collection of open sets in $M$ such that for each $V_j$ there exists a $U_{\nu(j)}$ such that $V_j\Subset U_{\nu(j)}$ and $\bigcup V_j= M$.

For a fixed $j$, the set $L^{-1}(L(V_j))$ is a union of $\mu_j$ connected components which are relatively compact in some open sets $U_{k_1},\ldots, U_{k_{\mu_j}}$; therefore, there exist $\mu_j$ linear maps $f_j^1,\ldots, f_j^{\mu_j}$ such that for each connected component there is one map which separates it from the others, that is, a map which has different values on it and on the union of the others.

Now, consider the map $F=(l_1,\ldots, l_{hp},f_1^1,\ldots, f_K^{\mu_K})$. By the first part of the construction, $F$ has an injective differential on $M$; by the second part, it is globally injective on $M$. Therefore $F$ is a holomorphic embedding of $M$ into $\C^{n}$, where $n=hp+\mu_1+\ldots+\mu_K$, realized with a complex linear map. $\Box$

\bigskip

Let $M$ be a compact $\Ci^1$ submanifold of a reflexive complex Banach space $E$ with complex structure $J$, with $\dim_\R M=2p-1$. $M$ induces a metric current $[M]$ of dimension $2p-1$.

\begin{Prp} \label{prp_caratt_CR}The following are equivalent:
\begin{enumerate}
\item $\dim_\C (T_z M\cap (JT_z M))=p-1$ $\forall\ z\in M$;
\item $[M](\alpha)=0$ for every metric $(r,s)-$form $\alpha$ on $E$ with $r+s=2p-1$ and $|r-s|>1$;
\item $M$ is locally the graph of a CR-function: for every $z\in M$ there exists $U$ neighbourhood of $z$ in $E$ such that $M\cap U$ is the graph of a function $f:\widetilde{M}\to E'$, $\widetilde{M}$ a CR-submanifold of $\C^p$, $E'$ a closed subspace of $E$ such that $E=E'\oplus \C^p$ and $f$ a CR-function on $\widetilde{M}$.
\end{enumerate}
\end{Prp}
\noindent{\bf Proof:} \emph{1) $\Longrightarrow$ 2) } 
Let $\alpha=(f,g_1,\ldots, g_r, h_1,\ldots, h_s)$ and let $i:M\to E$ an embedding whose differential is complex linear when restricted to (the preimage of) $T_zM\cap JT_zM$; then $[M]=i_\sharp T$, with $T\in M_{2p-1}(M)$. By the comparison theorem for manifolds, $T$ is induced by a classical current $T'$ on $M$; but then, 
$$[M](\alpha)=T(f\circ i, g_1\circ i,\ldots, h_s\circ i)=T'(f\circ i d(g_1\circ i)\wedge \ldots\wedge d(h_s\circ i))\;.$$
The functions $g_j\circ i$ have complex linear differentials on $i_*^{-1}(TM\cap JTM)$, therefore if there are more than $p$ of them, their wedge product will vanish; the same holds for the differentials of the functions $\overline{h}_j\circ i$. So $[M](\alpha)=0$ if $|r-s|>1$.

\medskip

\emph{2) $\Longrightarrow$ 1) } Let $\rho:E\to\C^N$ be a finite-dimensional embedding for $M$, which is holomorphic on $E$. If 
$$\dim_\C T_{\rho(z)}\rho(M)\cap J_{\C^N}T_{\rho(z)}\rho(M)=\dim_\C T_zM\cap JT_zM<p-1$$
then there exists a metric $(r,s)-$form $\beta$ on $\C^N$ with $r+s=2p-1$ and $|r-s|>1$ such that
$$\int_{\rho(M)}\beta\neq0$$
so
$$\int_M\rho^*\beta\neq0$$
and $\rho^*\beta$ is a $(r,s)-$form on $E$ with $|r-s|>1$.

\medskip

\emph{3) $\Longrightarrow$ 1) } Let $f:\widetilde{M}\to E'$ be the given CR-function; we define $G:\widetilde{M}\to \widetilde{M}\times E'$ by $G(p)=(p,f(p))$.

Let $F:\C^p\oplus E'\to E$ be the given isomorphism; then $(F\circ G)_*T_p\widetilde{M}=T_{F(p,f(p))} M$ and, since $T_{p}\widetilde{M}$ contains a complex subspace of complex dimension $p-1$, so does the tangent space of $M$.

\medskip

\emph{1) $\Longrightarrow$ 3) } Let us fix $z\in M$ and let $H_z$ be the complex subspace of $T_zE$ of (complex) dimension $p$ containing $T_zM$. By reflexivity, $E=E^*$,  so we have a splitting of $E=H_z\oplus E'$ for some closed subspace $E'$. By construction, $\pi:E\to H_z$ is a local embedding when restricted to a neighbourhood $U$ of $z$ in $M$, because it has a maximum rank differential at $z$.

Let $\widetilde{M}$ be the image of $U$ trough $\pi$; we have the function $f:\widetilde{M}\to E'$ defined by $(p,f(p))\in U\cap \pi^{-1}(p)$. By construction,  $f_*\vert_{T_p\widetilde{M}\cap JT_p\widetilde{M}}$ is $\C-$linear, so $f$ is a CR-function and $U$ is its graph. $\Box$

\bigskip

\begin{Def}Let $S$ be a $(2p-1)-$current with compact support in a complex manifold $X$; we say that $S$ is \emph{maximally complex} if $M_{r,s}=0$ for $|r-s|>1$.\end{Def}

\begin{Rem} $S_{r,s}$ in general won't be a metric current (see subsection 1.3 for an example). Nevertheless, the above definition makes sense for any current $S$, as we only require that the functional $M_{r,s}$ be zero for values $(r,s)$ such that $|r-s|>1$.\end{Rem}

\begin{Prp}\label{prp_slice_maxcmp}Let $M$ be a $(2p-1)-$current with compact support in $X$, $F:X\to Y$ a Lipschitz holomorphic map. Suppose that $M$ is maximally complex, then the same is true for $F_\sharp M$ and, if $p>\dim_\C Y$, for $\langle M, F, \zeta\rangle$, given that $M$ is flat and slices exist.\end{Prp}
\noindent{\bf Proof: }
Obviously, we have $(F_\sharp M)_{r,s}=F_\sharp(M_{r,s})$ (this is an equality between metric functionals only, not metric currents).

Moreover, if $\dim_\C Y<p$ and if $\langle M, F, \zeta\rangle$ exists for some $\zeta\in Y$, let $\{\rho_{\epsilon, \zeta}\}$ be a family of smooth approximations of $\delta_\zeta$. Then locally (with $\supp f$ contained in a manifold chart for $Y$)
$$\langle M, F, \zeta\rangle(f,\eta)=\lim_{\epsilon\to 0}M(f\rho_{\epsilon,\zeta}, F, \overline{F}, \eta)\;.$$
So, if $M_{r,s}=0$ for $|r-s|>1$ then also $\langle M, F, \zeta\rangle_{r-m, s-m}=0$ for $1<|r-s|=|(r-m)-(s-m)|$, with $m=\dim_\C Y$. $\Box$

\bigskip

\begin{Def}A \emph{MC-cycle} in a complex Banach space $E$ is a maximally complex $(2p-1)-$dimensional closed metric current, with compact support.\end{Def}

\begin{Rem}The definition is meaningless for $p=1$; the notion of \emph{moment condition} which substitutes the maximal complexity for $1-$dimensional currents cannot be given that easily in a Banach space and it turns out to be not automatically satisfied by a maximally complex current of higher dimension. The philosophical reason is the greater distance, in Banach spaces, between local and global aspects.\end{Rem}

The following Theorem follows easily from Proposition \ref{prp_slice_maxcmp} and from the slicing properties of rectifiable currents.

\begin{Teo} Let $M$ be a rectifiable MC-cycle of dimension $(2p-1)$ in a Banach space $E$ and consider a Lipschitz holomorphic map $F: E\to \C^m$. Then
\begin{enumerate}
\item $F_\sharp M$ is a rectifiable MC-cylce of dimension $(2p-1)$ in $\C^m$;
\item if $m<p-1$, $\langle M,F,\zeta\rangle$ is a rectifiable MC-cycle of dimension $2(p-m)-1$ in $E$.
\end{enumerate}
\end{Teo}

\begin{Rem}By Theorem \ref{teo_ext_slice}, the slice $\langle M, F, \zeta\rangle$ exists rectifiable for almost every $\zeta \in \C^m$.\end{Rem}

\begin{Teo} Let $M$ be a MC-cycle of dimension $(2p-1)$ in $E$. Then, for every linear projection $\pi:E\to \C^p$ and every $\phi\in\Ol(\overline{\supp M})$, we have
$$\debar[\pi_\sharp(\phi M)]^{0,1}=0\;.$$
Moreover, there is a unique integrable compactly supported function $c_\phi$ in $\C^p$ such that
$$\debar c_\phi=[\pi_\sharp(\phi M)]^{0,1}$$
and such a function can be obtained by convolution with the Cauchy kernel or the Bochner-Martinelli kernel.\end{Teo}
\noindent{\bf Proof: } We know that $dM=0$; since $M$ is maximally complex, we have
$$M=M_{p,p-1}+M_{p-1,p}$$
so
$$0=dM=(dM)_{p-2,p}+(dM)_{p-1,p-1}+(dM)_{p,p-2}\;.$$
In particular, this means that $(dM)_{p,p-2}=0$. Therefore
$$\debar[\pi_\sharp(\phi M)]^{0,1}=\debar[\pi_\sharp(\phi M)]_{p,p-1}=[d\pi_\sharp(\phi M)]_{p,p-2}=
[\pi_\sharp(d\phi M)]_{p,p-2}$$
$$=[\pi_\sharp(\phi dM)]_{p,p-2}+[\pi_\sharp(M\llcorner(1,\phi))]_{p,p-2}$$
but $M\llcorner(1,\phi)$ has non-vanishing $(r,s)-$components only for $(r,s)=(p-1,p-1)$ or $(r,s)=(p-2,p-1)$, so $[\pi_\sharp(M\llcorner(1,\phi))]_{p,p-2}=0$. Then
$$\debar[\pi_\sharp(\phi M)]^{0,1}=[\pi_\sharp(\phi dM)]_{p,p-2}=\pi_\sharp(\phi (dM)_{p,p-2})=0\;.$$
We note that $\pi_\sharp(\phi M)$ is a metric current in $\C^p$, therefore it is also a classical one, consequently its component of bidegree $(0,1)$ is a classical current as well and by the previous computation is $\debar-$closed. By a standard convolution-contraction with either the Cauchy kernel or the Bochner-Martinelli kernel, we can fin a compactly supported integrable function $c_\phi$ as requested. $\Box$

\medskip

\begin{Teo}\label{teo_HL1}Let $M$ be a compact, oriented $(2p-1)-$manifold (without boundary) of class $\Ci^2$ embedded in $H$, and suppose that there exists an orthogonal decomposition $H=\C^p\oplus H'$ such that the projection $\pi:H\to\C^p$, when restricted to $M$, is an immersion with transverse self-intersections. Then, if $M$ is an $MC-$cycle, there exists a unique holomorphic $p-$chain $T$ in $H\setminus M$ with $\supp T\Subset H$ and finite mass, such that $dT=[M]$ in $H$.\end{Teo}
\noindent{\bf Proof: } Let $\mathfrak{m}=\pi(M)\subset\C^p$; for every $\lambda\in H'\setminus\{0\}$, we define $\pi^\lambda(z)=(\pi(z), \langle z,\lambda\rangle)\in\C^{p+1}$.

By the previous results, $M^{\lambda}$ satisfies the same hypotheses in $\C^{p+1}$, therefore by \cite[Theorem 6.1]{harvey1} we can solve the problem for $M^{\lambda}=\pi^{\lambda}(M)$, finding a holomorphic $p-$chain $T^{\lambda}$ in $\C^{p+1}\setminus M$ with the required properties. Following the proof of Theorem 6.1 in \cite{harvey1}, we write
$$\C^p\setminus\mathfrak{m}=U_0\cup U_1\cup\ldots\cup U_k$$
where the $U_j$ are connected components and $U_0$ is unbounded; $T^\lambda$ is locally on each $U_j$ union of graphs of holomorphic functions
$$F_j^{\lambda, h}:U_j\to\C\qquad h=1,\ldots, n_{\lambda,j}\;.$$
Given another $\lambda'\in H'\setminus\{0\}$, we can consider the $p-$chain $T^{\lambda'}$, which will be given by holomorphic functions
$$F_j^{\lambda', h}:U_j\to\C\qquad h=1,\ldots, n_{\lambda',j}\;;$$
however, we can also consider, in $\C^{p+2}$, the manifold $M^{\lambda,\lambda'}$ and the associated solution $T^{\lambda,\lambda'}$; denoting by $p$ and $p'$ the restrictions of $\pi^\lambda$ and $\pi^{\lambda'}$ to $\C^{p+2}$, we have
$$p_* T^{\lambda,\lambda'}=T^{\lambda}\qquad p'_*T^{\lambda,\lambda'}=T^{\lambda'}\;.$$
Since the differentials of $p,\ p'$ are of rank $2p-1$ on $M^{\lambda,\lambda'}$ and because $p$ and $p'$ are holomorphic, their differentials are at least of rank $2p$ on $M^{\lambda,\lambda'}$; this means that they are of rank $2p$ in a neighborhood of $M^{\lambda,\lambda'}$ in $M^{\lambda,\lambda'}\cup\supp T^{\lambda,\lambda'}$ (which is locally a $\Ci^2$ manifold with boundary by Lemma 6.8 in \cite{harvey1}, ), therefore $n_{\lambda,j}=n_{\lambda,\lambda',j}=n_{\lambda',j}$ for every $j$ and every $\lambda,\lambda'\in H'\setminus\{0\}$.

Let $\{\lambda_i\}_{i\in I}$ be an orthonormal basis for $H'$ and consider the holomorphic functions
$$F_j^{\lambda_i, h}:U_j\to \C\qquad j=1,\ldots, k,\ \ h=1,\ldots, n_j,\ \ i\in I$$
and define
$$F_j^h=\sum_{i\in I}\lambda_i F_j^{\lambda_i, h}\;.$$

\noindent{\emph{The function $F_j^h$ is well defined. }} For any finite subset of indices $J\subset I$, we can consider the projection $$p_J:H\to\C^p\oplus\mathrm{Span}\{\lambda_i\}_{i\in J}$$
and the pushforward $[M]_J=(p_J)_\sharp [M]$; the functions $\{F^{\lambda_i, h}_j\}_{i\in J}$ give a solution for the finite-dimensional problem with datum $[M]_J$, therefore $S_{J,j,h}=\sum_{i\in J}F^{\lambda_i, h}_j\lambda_i$ is a holomorphic function with values in a finite-dimensional vector space, such that
$$|S_{J,j,h}(z)|\leq R$$
where $R$ is such that $\supp [M]_J\subset \C^p\times B(0,R)$, $B(z,r)$ being the ball with center $z$ and radius $r$ in $\mathrm{Span}\{\lambda_i\}_{i\in J}$.

Now, let us take $I=\N$ and fix $\epsilon>0$. By compactness, we can find $I'\subset I$ finite and set
$$V_\epsilon=\C^p\oplus\mathrm{Span}\{\lambda_i\}_{i\in I'}$$
so that $d(M,V_\epsilon)<\epsilon$; let $H'_\epsilon$ be the topological complement of $V_\epsilon$ in $H$, then the projection of $M$ on $H'_\epsilon$ lies in a ball of radius $\epsilon$ around $0$. Now, for any finite subset $J\subset I$ such that $\min J>\max I'$, we have that 
$$|S_{J,j,h}(z)|\leq \epsilon\;,$$
showing that the sequence of maps from $U_j$ to $H'$
$$\left\{\sum_{i=0}^m F^{\lambda_i,h}_j(z)\lambda_i\right\}_{m\in I}$$
is a Cauchy sequence with respect to the supremum norm on $U_j$. Therefore the limit $F^h_j(z)$ is well defined and continuous on the closure of $U_j$, because every element of the sequence is.

%
%
%

\medskip

\noindent{\emph{The function $F_j^h$ is holomorphic. } } Indeed, for any $\lambda\in H'$, we write 
$$\lambda=\sum_{i\in I}\alpha_i \lambda_i$$
and
$$\langle F_j^h(z), \lambda\rangle=\sum_{i\in I}\alpha_i F_j^{\lambda_i, h}(z)\;.$$
We now observe that
$$\left|\sum_{i\in I}\alpha_i F_j^{\lambda_i, h}(z)\right|\leq\sqrt{\sum_{i\in I}|F_j^{\lambda_i,h}|^2}\sqrt{\sum_{i\in I}|\alpha_i|^2}\leq\|\lambda\|_{H'}\sqrt{\sum_{i\in I}|F_j^{\lambda_i, h}(z)|^2_{\infty, U_j}}<\|\lambda\|_{H'}\||F^j_h|\|_{\infty, U_j}$$
which is finite, and this implies that the sequence of holomorphic functions
$$\left\{\sum_{i=0}^m\alpha_iF_j^{\lambda_i,h}\right\}_{m\in I}$$
converges uniformly on $U_j$. The limit is then holomorphic, so $F_j^h$ is holomorphic.

\medskip

\noindent{\emph{The function $F_j^h$ extends $\Ci^1$ to the boundary. } } By \cite{harvey1}, there exist sets $A\subset\mathfrak{m}$ and $A_i\subset M^{\lambda_i}$ with $\pi(A_i)=A$, which are $\H^{2p-1}-$negligible and such that outside them we have $\Ci^1$ regularity for $\supp T^{\lambda_i}\cup M^{\lambda_i}$ and for the functions $F_j^{\lambda_i,h}$. Let us consider $p\in \mathfrak{m}\cap \overline{U_j}\setminus A$; for each $i\in I$, one of the following two cases can occur:
\begin{enumerate}
\item $F_j^{\lambda_i, h}(p)\not\in M^{\lambda_i}$,
\item $F_j^{\lambda_i,h}(p)\in M^{\lambda_i}$.
\end{enumerate}
In the former, $F_j^{\lambda_i,h}$ extends holomorphically through $p$, whereas in the latter we can find a relatively compact neighborhood $V$ of $p$ in $\mathfrak{m}$ such that $F_{j}^{\lambda_i,h}$ coincides on $V$ with some $CR$ function $f:V\to M^{\lambda_i}$. In both cases, $F_j^{\lambda_i,h}$ is of class $\Ci^1$ near $p$. Let $U$ be an open set with $\Ci^1$ boundary in $U_j$ such that $bU_j\cap bU= V$.

The restrictions of the derivatives of $F^{h}_j$ to $bU_j$ are continuous, when we derive in a direction tangent to $TbU$; however, by the Cauchy-Riemann equations, we can control the normal derivative with the tangential ones, therefore also the normal derivative of $F^{h}_j$ is a continuous function when restricted to $bU$.

We note that from this follows that the image of $bU$ through one of these maps is a compact set in $H'$ and we can replicate the previous argument, obtaining that the sequence
$$\left\{\sum_{i=0}^m \frac{\de}{\de z_s}F^{\lambda_i,h}_j(z)\lambda_i\right\}_{m\in I}$$
is a Cauchy sequence with respect to the supremum norm on $U$.

Therefore, the limit is continous on the closure of $U$, thus implying that
$$\left\|\left|\frac{\de}{\de z_s}F^h_j\right|\right\|_{\infty, U}<+\infty\;.$$
Moreover, on $bU\cap bU_j=V$, $F^h_j$ coincides with $f$ and we can cover $\H^{2p-1}-$almost all of $bU_j$ with open sets where $F^h_j$ coincides with some CR-functions realizing $M$ as a graph. Therefore, as $M$ is a compact $\Ci^1$ manifold,
$$\left\|\left|\frac{\de}{\de z_s}F^h_j\right|\right\|_{\infty, bU_j}<+\infty$$
hence
$$\left\|\left|\frac{\de}{\de z_s}F^h_j\right|\right\|_{\infty, U_j}<+\infty\;.$$

\medskip

\noindent{\emph{The current of integration on the graph of $F_j^h$ has finite mass. } } By the previous paragraph, there exists a constant $C_{h,j}$ such that
$$|\nabla F^h_j(z)|^2=\sum_{i\in I}|\nabla F^{\lambda_i, h}_j(z)|^2\leq C_{h,j}\qquad\textrm{ for every }z\in U_j\;.$$
It is easy to show that there exists a polinomial $g_p(X)$ such that
$$\sum_i |a_i|\leq S<+\infty\Longrightarrow \sum_{|J|=p}\prod_{i\in J}|a_i|\leq g_p(S)<+\infty\;.$$
Therefore
$$\sum_{|J|=p}\prod_{i\in J}|\nabla F^{\lambda_i, h}_{j}(z)|\leq g_p(C_{h,j})<+\infty\qquad\textrm{ for every }z\in U_j\;.$$
We consider the $(p,p)-$form
$$\eta^h_j(z)=\sum_{|J|=p}\bigwedge_{i\in J}dF^{\lambda_i, h}_j(z)\wedge d\overline{F}^{\lambda_i, h}_j(z)$$
which is well-defined by the previous estimates and note that
$$\|\eta^h_j\|_{\infty, U_j}\leq g_p(C_{h,j})\;.$$
Let $\{w_i\}_{i\in I}$ be coordinates for the basis $\{\lambda_i\}_{i\in I}$, i.e. $w_i(v)=\langle v,\lambda_i\rangle$ for $v\in H'$, and denote by $T_{h,j}$ the (alleged) current of integration on the graph of $F^h_j$. Then
$$T_{h,j}(1, w_{i_1}, \overline{w}_{i_1},\ldots, w_{i_p},\overline{w}_{i,p})=\int\limits_{U_j}dF^{\lambda_{i_1}, h}_j(z)\wedge d\overline{F}^{\lambda_{i_1}, h}_j(z)\wedge\ldots\wedge dF^{\lambda_{i_p}, h}_j(z)\wedge d\overline{F}^{\lambda_{i_p}, h}_j(z)\;.$$
Therefore, by Proposition \ref{prp_vol_in_H}, we have
$$M(T_{h,j})\leq \mathcal{L}^{2p}(U_j)\sum_{p'=0}^pg_{p'}(C_{h,j})<+\infty\;.$$
We have to sum all the values from $0$ to $p'$ because we apply the formula of Proposition \ref{prp_vol_in_H} in $H$ and not in $H'$, so we have to consider also the $p-$tuples of coordinates coming in part from $\C^p$ and in part from $H'$.

\medskip

As the $F_j^h$ are a finite number of functions, we can consider the metric functional of integration on their graphs and denote it by $T$. $T$ is a holomorphic $p-$chain in $H\setminus M$, it has finite mass and its support is contained in a product of discs, therefore it is relatively compact in $H$. Moreover, for $\H^{2p-1}-$almost every point in $M$ there is a neighborhood where $\supp T\cup M$ is a $\Ci^1$ manifold.

This implies that $T$ is a metric rectifiable $(p,p)-$current in $H$. We note that for any finite-dimensional projection $p:H\to\C^m$, we have that $d(p_\sharp T)=p_\sharp[M]$; it is an easy application of Theorem \ref{teo_ext} to show that this implies $dT=[M]$.  Finally, it is not difficult to see that the map $x\mapsto (x,F_j^h(x))$ is proper into $H\setminus M$, which is an hilbertian manifold, hence by Theorem \ref{teo_hilb_img} its image is a finite dimensional complex space in $H\setminus M$.  
$\Box$

\bigskip

\begin{Rem}Suppose we are given a family $M_s$ of $MC-$cycles each satisfying the hypothesis of Theorem \ref{teo_HL1}, depending on some parameter $s\in U\subset \C$ in a $\Ci^1$ way. Locally in $s$, we can assume that the various manifolds $M_s$ project, through the maps $\pi_s$ given by hypothesis, to the same immersed manifold $\mathfrak{m}\subset\C^p$.

Therefore, the functions $F^{\lambda, h}_{j,s}$ constructed during the proof of Theorem \ref{teo_HL1} depend on $s$ in a $\Ci^1$ way. This implies that $F^{h}_{j,s}$ vary continuously in $s$.

Hence, the map associating to $s$ the solution to $dT_s=[M_s]$ is continuous in $s$.\end{Rem}

\medskip

The compactness of $M$ is needed only to ensure that there is a finite number of connected components in $\C^p\setminus\mathfrak{m}$. Therefore, we also have the following result.

\begin{Teo}\label{teo_HL2}Let $M$ be a bounded, oriented $(2p-1)-$manifold (without boundary) of class $\Ci^2$ embedded in $H$, with finite volume.  Assume that there exists an orthogonal decomposition $H=\C^p\oplus H'$ such that the projection $\pi:H\to\C^p$, when restricted to $M$, is a closed immersion with transverse self-intersections. Then, if $M$ is an $MC-$cycle, there exists a unique holomorphic $p-$chain $T$ in $H\setminus M$ with $\supp T\Subset H$ and finite mass, such that $dT=[M]$ in $H$.\end{Teo}
\noindent{\bf Proof: }As $\pi$ is supposed to be closed, $\pi(M)$ is a closed and bounded subset of $\C^p$, therefore compact. By the finiteness of volume, we know that $[M]$ is a well-defined metric current and we can proceed with the same proof as before. $\Box$

\bibliography{bibsing}{}

\begin{thebibliography}{10}

\bibitem{ambrosio1}
{\sc L.~Ambrosio and B.~Kirchheim}, {\em Currents in metric spaces}, Acta
  Math., 185 (2000), pp.~1--80.

\bibitem{ambrosio2}
\leavevmode\vrule height 2pt depth -1.6pt width 23pt, {\em Rectifiable sets in
  metric and {B}anach spaces}, Math. Ann., 318 (2000), pp.~527--555.

\bibitem{ambrosio3}
{\sc L.~Ambrosio and T.~Schmidt}, {\em Compactness results for normal currents
  and the {P}lateau problem in dual {B}anach spaces}.
\newblock http://cvgmt.sns.it/paper/1762/, 2012.

\bibitem{aurich1}
{\sc V.~Aurich}, {\em Local analytic geometry in {B}anach spaces}, in Complex
  analysis, functional analysis and approximation theory ({C}ampinas, 1984),
  vol.~125 of North-Holland Math. Stud., North-Holland, Amsterdam, 1986,
  pp.~1--23.

\bibitem{douady1}
{\sc A.~Douady}, {\em Le probl\`eme des modules pour les sous-espaces
  analytiques compacts d'un espace analytique donn\'e}, Ann. Inst. Fourier
  (Grenoble), 16 (1966), pp.~1--95.

\bibitem{harvey1}
{\sc F.~R. Harvey and H.~B. Lawson, Jr.}, {\em On boundaries of complex
  analytic varieties. {I}}, Ann. of Math. (2), 102 (1975), pp.~223--290.

\bibitem{harvey3}
\leavevmode\vrule height 2pt depth -1.6pt width 23pt, {\em On boundaries of
  complex analytic varieties. {II}}, Ann. of Math. (2), 106 (1977),
  pp.~213--238.

\bibitem{harvey2}
{\sc R.~Harvey and B.~Shiffman}, {\em A characterization of holomorphic
  chains}, Ann. of Math. (2), 99 (1974), pp.~553--587.

\bibitem{king1}
{\sc J.~R. King}, {\em The currents defined by analytic varieties}, Acta Math.,
  127 (1971), pp.~185--220.

\bibitem{mongodi1}
{\sc S.~Mongodi}, {\em Some applications of metric currents to complex
  analysis}.
\newblock to appear in \emph{Manuscripta Mathematica}, 2012.

\bibitem{mujica1}
{\sc J.~Mujica}, {\em Complex analysis in {B}anach spaces}, vol.~120 of
  North-Holland Mathematics Studies, North-Holland Publishing Co., Amsterdam,
  1986.
\newblock Holomorphic functions and domains of holomorphy in finite and
  infinite dimensions, Notas de Matem{\'a}tica [Mathematical Notes], 107.

\bibitem{noverraz1}
{\sc P.~Noverraz}, {\em Pseudo-convexit\'e, convexit\'e polynomiale et domaines
  d'holomorphie en dimension infinie}, North-Holland Publishing Co., Amsterdam,
  1973.
\newblock North-Holland Mathematics Studies, No. 3. Notas de Matem{\'a}tica
  (48).

\bibitem{ruget1}
{\sc G.~Ruget}, {\em \`{A} propos des cycles analytiques de dimension infinie},
  Invent. Math., 8 (1969), pp.~267--312.

\bibitem{shiffman1}
{\sc B.~Shiffman}, {\em On the removal of singularities of analytic sets},
  Michigan Math. J., 15 (1968), pp.~111--120.

\end{thebibliography}
\bibliographystyle{siam}

\end{document}